\input amstex
\documentstyle{amsppt}
\magnification=1200
\hsize=13.8cm
\vsize=19.0 cm
\catcode`\@=11
\def\NoLogo{\let\logo@\empty}
\catcode`\@=\active \NoLogo

\def\qeddemo{\qed\enddemo}

%\long\def\MM#1{{\smallskip\noindent\tt ** #1\smallskip}}
\def\today{\ifcase\month\or  %from ltr.sty
 January\or February\or March\or April\or May\or June\or
 July\or August\or September\or October\or November\or December\fi
 \space\number\day, \number\year}
\def\heat{\left(\frac{\p}{\p t}-\Delta\right)}
\def\covar{\nabla}
\def\Ric{Rc}
\def\Rscalar{R}
\def\a{\alpha}
\def\p{\partial}
\def\pa{\partial}

\def\R{\Bbb R}
\def\D{\Delta}

\def\ctL{{\Cal L_+}}
\def\cW{{\Cal W}}
\def\cF{\Cal F}
\def\cN{\Cal N}
\def\ctN{\Cal N_+}
\def\ctF{\Cal F_+}
\def\ctW{\Cal W_+}
\def\Vol{\operatorname{vol}}

%formulas
\def\RF{1.1}
\def\WPdef{1.2}
\def\Ueq{1.3}
\def\WPmono{1.4}
\def\WPP{1.5}
\def\VPP{1.6}
\def\Vmono{1.7}
\def\WPestimate{1.8}

%theorems
\def\WPmonoprop{1.1}
\def\Harnackprop{1.2}
\def\Blowdownthm{1.3}
\def\WPboundedprop{1.4}
\def\WPestimateprop{1.5}
\def\Hauptvermutung{1.6}
\def\Munuprop{1.7}

\def\PosCurvOp{2.1}
\def\ExpanderEll{2.2}
\def\ThetaLower{2.3}

\leftheadtext{Feldman, Ilmanen, and Ni}
\rightheadtext{Expander Entropy}

\topmatter
\title{Entropy and reduced distance for Ricci expanders}
\endtitle
\author{Michael Feldman \footnotemark, Tom Ilmanen
\footnotemark  and Lei Ni\footnotemark }\endauthor
\footnotetext"$^{1}$"{Department of Mathematics, University of Wisconsin,
Madison WI 53706. Research partially supported by NSF grant DMS-0200644.}
\footnotetext"$^{2}$"{Departement Mathematik, ETH Zentrum, 8092 Z\"urich,
Switzerland. Research partially supported by
Schweizerische Nationalfonds grant 21-66743.01.}
\footnotetext"$^3$"{Department of Mathematics, 
University of California, San Diego, La
Jolla CA 92093. Research partially supported by NSF grant
DMS-0328624 and an Alfred P. Sloan Fellowship.}

\address
%Departement Mathematik, ETH, Z\"urich
\endaddress
\email{}
\endemail

\address
%Department of Mathematics, University of Wisconsin
\endaddress
\email{}
\endemail

\address
%Department of Mathematics, University of California, San Diego, La
%Jolla, CA 92093
\endaddress
\email{ lni\@math.ucsd.edu}
\endemail

\address
Department of Mathematics,
\endaddress
\email{}
\endemail

\affil 
\endaffil

\date\today\enddate

\abstract 
Perelman has discovered two integral quantities, the
shrinker entropy $\cW$ and the (backward) reduced volume, that are
monotone under the Ricci flow $\pa g_{ij}/\pa t=-2R_{ij}$ and constant
on shrinking solitons.  Tweaking some signs, we find similar formulae
corresponding to the expanding case. The {\it expanding entropy}
$\ctW$ is monotone on any compact Ricci flow and constant precisely on
expanders; as in Perelman, it follows from a differential inequality
for a Harnack-like quantity for the conjugate heat equation, and leads
to functionals $\mu_+$ and $\nu_+$.  The {\it forward reduced volume}
$\theta_+$ is monotone in general and constant exactly on
expanders.

A natural conjecture asserts that $g(t)/t$ converges as $t\to\infty$
to a negative Einstein manifold in some weak sense (in particular
ignoring collapsing parts).  If the limit is known a-priori to be
smooth and compact, this statement follows easily from any monotone
quantity that is constant on expanders; these include
$\Vol(g)/t^{n/2}$ (Hamilton) and $\bar\lambda$ (Perelman), as well as
our new quantities.  In general, we show that if $\Vol(g)$ grows like
$t^{n/2}$ (maximal volume growth) then $\ctW$, $\theta_+$ and
$\bar\lambda$ remain bounded (in their appropriate ways) for all
time. We attempt a sharp formulation of the conjecture.
\endabstract
\endtopmatter

\document
\vskip .2cm

Small, large and distant parts of a Ricci flow are known to be
modelled by various kinds of Ricci solitons: steady,
shrinking, and expanding. Perelman has discovered monotone quantities
correponding to the first two of these.  The functional $\cF$ is monotone
on the Ricci flow and constant precisely on steadies \cite{P1,
\S1}. The shrinking, or localizing, entropy $\cW$ is monotone in
general and constant precisely on shrinking solitons \cite{P1, \S3}.
Closely related is the notion of the (backward) reduced volume
\cite{P1, \S7}.
The latter two show that developing singularities and
ancient histories are modelled on shrinkers, when an appropriate
blowup or blowdown is taken \cite {P1, S}.

In this note, by tweaking some signs, we observe similar monotonicity
formulae for the expanding case.  In \S1, analogous to \cite{P1, \S1},
we define an {\it expanding}, or {\it delocalizing entropy} $\ctW$,
which is monotone in general and constant precisely on expanders. It
follows from a differential inequality for the conjugate heat equation
similar to \cite{P1, \S9}; see also \cite{LY}. We also construct
functionals $\mu_+(g,\sigma)$ and $\nu_+(g)$ from $\ctW$.

In \S2, analogous to \cite{P1, \S7}, we present a {\it forward reduced
distance} $\ell_+$ and its corresponding {forward reduced volume}, which
is again monotone in general, and constant exactly on
expanders. 

Our formulae may be compared to those for $\Vol(g)/t^{n/2}$ (Hamilton)
and of $\bar\lambda(g):=\Vol(g)^{2/n}\lambda(g)$ (Perelman). All of
these monotonicities have implications for the convergence of the
rescaled Ricci flow as $t\to\infty$. If $g(t)$ is a compact Ricci flow
that exists for all $t>0$, and the volume $\Vol(g)$ grows like
$t^{n/2}$ (maximal volume growth), we observe in Theorem
\WPboundedprop\ that $\ctW$ remains bounded. A similar result holds
for Perelman's scaled $\lambda$-functional, and for the forward reduced
volume (Proposition \ThetaLower). This suggests that
(ignoring zero-volume collapsing parts) the long-term limit of
$g(t)/t$ is a negative Einstein manifold, possibly with singularities;
see Conjecture \Hauptvermutung.

Our results may be set against recently discovered entropy formulae
and pointwise differential inequalities in other settings: for the
standard heat equation on a manifold of nonnegative Ricci curvature
\cite{N1}, and for the forward conjugate heat equation $u_t=\Delta
u+Ru$ on a 2-dimensional Ricci flow \cite{N2} and (matrix Harnack!) on
a K\"ahler-Ricci flow \cite{N3}.\footnote{It is conjugate to the
backward heat equation $u_t=-\Delta u$.}

It is also instructive to compare the Ricci flow to the mean curvature
flow in Euclidean space.  The monotonicity formulae for $\cF$, $\cW$,
$\ctW$, and for the corresponding reduced volumes, are analogous to
the following three monotone energies on a mean curvature flow $M_t$
in $\R^N$:
$$
\gather
\frac{d}{dt}\int_{M_t} 
\frac{e^{-|x-x_0|^2/4(t_0-t)}}{(-4\pi t)^{n/2}}\le0,\qquad t<0,\\
\frac{d}{dt}\int_{M_t} e^{\langle x,x_0\rangle
-|x_0|^2t}\le0,\qquad t\in\R,\qquad\qquad
\frac{d}{dt}\int_{M_t} 
\frac{e^{|x-x_0|^2/4(t-t_0)}}{(4\pi t)^{n/2}}\le0,\qquad t>0.
\endgather
$$ 
The first is Huisken's well-known result \cite{Hu}
and is constant precisely on shrinkers.  The second picks out steadies
(translators) \cite{I1}, whereas the
third is stationary on expanders \cite{I4}.

In the noncompact case, for both flows, the energy of a steady or
expander is typically infinite, so in this case the monotonicity
formulae will not be useful in proving convergence to an asymptotic
soliton.\footnote{There are graphs over $\R^n$ (due to Ecker)
that evolve by mean curvature flow forever, but are
asymptotic neither to expanders nor translators.  There also appear to
be noncompact curves in $\R^2$ whose mean curvature flow exhibits
nontrivial periodicity, i.e.\ breathers \cite{BI, I4}.}
This is unfortunate, because noncompact expanders are 
needed for modelling the relaxation out of a previously formed
singularity of conic type; see for example \cite{ACI, I3, FIK}
and remarks after Proposition \ThetaLower.

In the compact case, a sharp disanalogy between the two flows
arises. Every compact mean curvature flow 
disappears in finite time, so mean curvature flow has no
compact steadies or expanders. 
On the other hand, in Ricci flow, negative curvature allows
a compact solution to last forever, even grow.
So compact steadies and expanders exist, but on sufferance:
they are always Einstein \cite{H3, \S2; CK Prop. 5.20}. 

On the other hand, shrinkers typically have finite energy whether
compact or noncompact. Furthermore, the soliton potential need not be
constant. (So a Ricci shrinker $g(t)$ need not be Einstein, and a mean
curvature shrinker $M_t$ need not be a minimal submanifold of a sphere.)
There is accordingly a lively variational calculus for them (see \cite
{Hu, I2, I4, W, CHI} and others).

\bigskip
\subheading{\S1 Expander entropy formula}

\bigskip
Let $(M, g(t))$ be a solution to the Ricci flow
$$
\frac{\p g_{ij}}{\p t} =-2R_{ij}. \tag\RF
$$
We call $g(t)$ a {\it (gradient) expander} if 
$$
g(t)=a(t)\phi^*_t(g_0),
$$
where $\phi_t$ is a family of diffeomorphisms whose velocity
is given by a gradient vector field, and $a(t)$ is any 
increasing function. Then $a(t)=t-T$ for some $T$ 
and $g$ solves the PDE
$$
Rc+\covar^2F+\frac{g}{2(t-T)}=0,\qquad t>T,
$$
on each time slice, for some function $F=F(\cdot,t)$. 

Now let $M$ be closed, $g$ a metric, $u>0$ a function,
and $\sigma>0$. Dual to the shrinker entropy $\cW$ \cite{P1, Section 3},
(which we call $\cW_-$) define the {\it expander entropy} by
$$
\split
\ctW(g,u,\sigma)&:=
\int_M \left[\sigma\left(\frac{|\covar u|^2}{u}+\Rscalar u\right)
+u\log u\right]\, dv +\frac{n}{2}\log(4\pi\sigma)+n,\\
&=\int_M \left[\sigma(|\covar f_+|^2+{\Rscalar})
-f_++n\right]u\, dv \\
\endsplit  
\tag\WPdef
$$
restricted to $u$ satisfying
$$
\int_M u\, dv=1.
$$
Here $f_+=f_+^{u,\sigma}$ is defined via 
$$
u=\frac{e^{-f_+}}{(4\pi\sigma)^{n/2}}
$$
Now take $u=u(x,t)$ to be a positive solution to the {\it
conjugate heat equation}
$$
\frac{\p u}{\p t}=-\D  u +{\Rscalar}u, \tag\Ueq
$$
satisfying the conserved condition $\int u\, dv=1$.
Note that $f_+=f_+^{u,t-T}$ solves 
$$
\frac{\p f_+}{\p t}=-\D  f_+ +|\covar  f_+|^2 -{\Rscalar}-\frac{n}{2
(t-T)}.
$$
The following monotonicity formula for $\ctW$ was found by varying
the signs in the shrinker entropy formula \cite{P1, \S3}.  A more
elegant derivation is given a few pages hence.

\proclaim{Theorem \WPmonoprop\ (Expander entropy monotonicity)} 
Under evolution equations (\RF), (\Ueq) we have
$$
\frac{\p}{\p t}\ctW(g(t),u(t),t-T)=\int_M
2(t-T)u\left|R_{ij}+\covar_i\covar_jf_++\frac{g_{ij}}{2(t-T)}\right|^2. 
\, dv, \tag\WPmono
$$
%where $f_+=f_+_{u,\sigma}$ is defined via 
%$$
%u={e^{-f_+}}/(4\pi\sigma)^{n/2}.
%$$
\endproclaim

\noindent
Note that the right hand side vanishes precisely on an expander
with birth time $T$.

Formula (\WPmono) is a consequence of the following differential
inequality, which is similar to the one in \cite{P1, \S9} (see also
\cite{LY}), but with sign changes.

\proclaim{Theorem \Harnackprop} Let $(M, g(t))$, $u$, $f_+$, and $T$
be as above. Define
$$
v_+:=\left[(t-T)\left(2\D f_+-|\covar 
f_+|^2+{\Rscalar}\right)-f_++n\right]u,\qquad t>T.
$$
Then
$$
\left(\frac{\p}{\p t}+\D -{\Rscalar}\right)v_+
=2(t-T)u\left|R_{ij}+\covar_i\covar_jf_++\frac{g_{ij}}{2(t-T)}\right|^2,\qquad t>T.
\tag \WPP
$$
\endproclaim

\demo{Proof} 
We may assume $T=0$.  Set $V:=2\D
f_+-|\covar  f_+|^2 +{\Rscalar}$ and compute
$$
\left(\frac{\p}{\p t} +\D\right)V=2|R_{ij}+\covar_i\covar_jf_+|^2+2\langle
\covar V, \covar f_+\rangle. \tag \VPP
$$
and
$$ 
\left(\frac{\p}{\p t} +\D\right)(t V-f_+ +n)=2t
|R_{ij}+\covar_i\covar_jf_+|^2 +2\langle\covar (t V-f_++n),\covar f_+\rangle
+2\Delta f_++2R+\frac{n}{2t}
$$
and do further calculation.
\qeddemo

\demo{Proof of Theorem \WPmonoprop}
Note that $\int v_+\, dv = \ctW$. To get the time derivative of $\ctW$,
integrate (\WPP) over $M$.
\qeddemo

\demo{Remark}
For the steady case, one can define $e^{-f}=u$,
$v_0:=Vu=(2\Delta f-|\covar f|^2+R)u$ and verify 
$$
\left(\frac{\p}{\p t}+\D -{\Rscalar}\right)v_0
=2u\left|R_{ij}+\covar_i\covar_jf\right|^2,
$$
which may be viewed as a ``neutral'' differential Harnack inequality.
It may be integrated to get the monotonicity of $\cF$ in \cite{P1,
\S1}. See \cite{CCLN}. 
\qeddemo

We next consider the implications for long term behavior 
of the Ricci flow. 
Let $M$ be compact and consider a Ricci flow $g(t)$ that exists for
all $t>b$. We construct a global solution of the conjugate heat equation.
Let $t^i\to\infty$, and let $u^i$, $b<t\le t^i$ be the solution of (\Ueq)
with final data
$$
u^i(x,t^i)=\frac{1}{V(t^i)},\qquad x\in M.
$$
Passing $t^i\to\infty$ we obtain
an immortal  solution $u>0$ of (\Ueq) with $\int_M u=1$ for all $t>b$.

Now let $\a_j\to\infty$ and define the {\it blowdown sequence}
$$
g_j(t):=\a_j^{-1}g(\a_jt),\qquad t> b/\a_j.
$$
The following theorem is classic. We give a new proof using $\ctW$.
The hypothesis will hold (subsequentially) if, for example, 
$|Rm|$, diameter, and injectivity radius of $g(t)/t$ 
are controlled for all time \cite{H2}.
A similar theorem has been proven for the shrinking case \cite{S}.

%The curvature bound is called a type III condition.

\proclaim{Theorem \Blowdownthm\ (Hamilton \cite{H3})} 
Let $M$ be compact. Suppose there are diffeomorphisms $\phi_j$ such that
$\phi^*_j(g_j(t))\to g_\infty(t)$ smoothly, $t>0$. 
Then $g_\infty(t)$ is an expanding soliton (hence negative Einstein).
\endproclaim

\demo{Proof} Let $u_j(t):=\a_j^{-n/2}u(\a_jt)$. Evidently we can
extract a smoothly convergent subsequence $\phi_j^*(u_j)\to
u_\infty>0$ satisfying (\Ueq) for all time. The convergence is smooth,
the manifold is compact, and $\ctW$ is monotone and scale
invariant. So $\ctW$ is constant on the limit flow, which is therefore
an expanding soliton.  \qeddemo

\demo{Remark} 
Note that an expanding breather can be expressed as a
blowdown of itself, after extending it iteratively to a maximal
interval of existence $t>b$. Applying the Theorem then 
yields an alternative proof that
there are no nontrivial compact expanding breathers \cite{P1, \S2}.
\qeddemo

Theorem \Blowdownthm\ may be proven using any monotone, scale invariant
quantity that is constant only on expanders.  Hamilton's original
method \cite{H3} uses the quantity $R_{min}(t)$ $V(t)^{2/n}$. Alternately
(Hamilton) one can use the {\it scaled volume}
$$
\tilde V(t):=\frac{V(t)}{t^{n/2}}.
$$
Here $V(t)$ is the volume of $g(t)$ and $\tilde V(t)$ 
is the volume of the rescaled metric $\tilde g(t):=g(t)/t$.
Perelman's proof \cite{P1, \S2} uses the {\it scaled
$\lambda$-functional}
$$
\bar\lambda(g):=V(g)^{2/n}\lambda(g),
$$ 
where $\lambda(g):=\inf\{\cF(g,u)|\int u\,dv=1\}$.
One can also use the forward reduced volume (see \S2).
In each case, the right-hand side of the monotonicity formula
vanishes identically only on expanders.

We now drop the hypothesis that the blowdowns converge smoothly to a
compact limit. Can we still get convergence to expanders?
The various monotonicity formulae will only be useful if the monotone quantity
remains bounded for all time. Accordingly, we embark on a study
of the long-term behavior of $\ctW$, $\tilde V$, and $\bar \lambda$.

We begin with the scaled volume. Recall from Hamilton \cite{H3} 
$$
\frac{\p R}{\p t}=\Delta R+2|Rc|^2,
$$
from which follows via $R^2\le n|Rc|^2$ and the maximum principle that
$$
R+\frac{n}{2t}\ge0,\qquad
\frac{d\tilde V}{dt}
=-\frac{1}{t^{n/2}}\int\left(R+\frac{n}{2t}\right)\,dv \le0.
\tag\Vmono
$$ 
So $\tilde V(t)$ is nonincreasing, and $\tilde
V_\infty:=\lim_{t\to\infty}\tilde V(t)$ exists.  \footnote{The 
proof of Theorem \Blowdownthm\ (Hamilton) then runs as follows: assume $g_j(t)$
converges smoothly to a compact limit $g_\infty(t)$.  Then $\tilde
V(g(t))\to\tilde V(g_\infty(t))$, which is therefore constant, so
$R_\infty\equiv-n/2t$, so $Rc_\infty$ is diagonal and equals
$-g_\infty/2t$. This shows that any such blowdown is 
negative Einstein. In particular, any compact expander
is negative Einstein.}

Next we motivate $\ctW$ as mentioned above. Analogously to 
\cite{P1, \S5} (see also \cite{N2}), define 
$$
\cN(g,u):=\int_M u\log u\, dv, 
\qquad
\cN_+(g,u,\sigma):=\cN(g,u)+\frac{n}{2}\log (4\pi
\sigma)+\frac{n}{2}.
$$ 
The quantity $\cN$ is the {\it Nash entropy}
\cite{Na}\footnote{Our sign convention for $\cN$ is the opposite of
Nash's.}. Note that $\cN_+$ is dual to the entropy 
used for shrinkers in \cite{P1, \S5}.
Now let $g(t)$ and $u(t)$ evolve as above for all $t>0$ 
and set $\sigma(t):=t$. Compute
$$
\frac{d\cN}{dt}=\cF:=\int_M \left(\frac{|\covar u|^2}{u}+Ru\right)\,dv,\qquad
\frac{d\cN_+}{dt}=\cF_+:=\cF+\frac{n}{2t},
$$
and
$$
\split
\frac{\p}{\p t}(t\cN_+)&= \cN_++t\cF_+\\
&=\int u\left[\log u+\frac{n}{2}\log(4\pi t)+\frac{n}{2}+\frac{n}{2}
+t\left(\frac{|\covar u|^2}{u}+Ru\right)\right]\,dv\\
&={\ctW},\\
\endsplit
$$
which motivates $\ctW$.

Next we recall \cite{P1, \S2}
$$
\frac{d}{dt}\cF\ge\frac{2}{n}\cF^2.
$$ 
%In fact, the $u$ that minimizes $\cF$ is unique and depends smoothly
%on $g$, so $\lambda(g(t))$ is smooth and satisfies
%$$
%\frac{d\lambda}{dt}\ge \frac{2}{n}\lambda^2.
%$$
By ODE comparison, if $g(t)$ exists for all $t>0$ we have
$$
-\frac{n}{2t}\le\cF(g,u)\le0,\qquad t>0.
$$
(Side remark: passing the starting time to $-\infty$, we notice that
a compact eternal solution has $\cF\equiv0$, hence is a steady soliton.)

Next we study $\cN$, $\ctW$, and $\bar\lambda$ as $t\to\infty$.
By Jensen's inequality, we have in general
$$
\int u\log u\,dv\ge \log\left(\frac{1}{V}\right)
$$
when $\int u\,dv=1$, with equality precisely when $u=1/V$.
So we have the lower bound
$$
\cN_+(g,u,t)
\ge -\log(V(t))+\frac{n}{2}\log(4\pi t)
+\frac{n}{2}.
$$
To get a similar upper bound, we use monotonicity and the 
construction of $u$ as asymptotically constant. 
From above we have 
$$
\frac{d}{dt}\ctN(g, u, t)=\ctF(g,u, t) =\cF(g,u)+\frac{n}{2t}\ge 0,
$$
so $\cN_+$ is monotone nondecreasing. This also applies also to $u^i$. 
So fixing $t$, we have for all large $t^i$,
$$
\ctN(g, u^i, t)\le\ctN(g, u^i,t^i)
=\log\left(\frac{1}{V(t^i)}\right)
+\frac{n}{2}\log(4\pi t^i)+\frac{n}{2}
$$
Taking $i\to\infty$, combining with the lower bound, and sending $t\to\infty$,
we obtain
$$
\lim_{t\to\infty}\ctN(g, u, t)=
-\log(\tilde{V}_{\infty})+\frac{n}{2}\left(\log4\pi+1\right).
$$
(This holds even if $\tilde V_\infty=0$.) If 
$\tilde V_\infty>0$, then since $\ctN$ is monotone,
there exist $t_j \to \infty$ such
that 
$$ 
\lim_{j\to \infty} t_j \ctF(g, u, t_j)=0.
$$
Since $\ctW=t\ctF+\ctN$, and since $\ctW$ is monotone, we obtain
$$
\lim_{t\to\infty}\ctW(g, u, t)=\lim_{t\to\infty}\ctN(g, u, t).
$$
(Again, it holds even if $\tilde V_\infty=0$.)

Next we turn to $\bar\lambda$. It is straightforward
to show that the minimizer in the definition of $\lambda$
depends smoothly on $g$, so $\lambda(g(t))$ is smooth,
and 
$$
\frac{d\lambda}{dt}\ge \frac2n\lambda^2,\qquad-\frac{n}{2}\le t\lambda\le0,
$$
since $g(t)$ exists for all $t>0$. 
From $t_j\ctF\to0$ and the lower bound for $\lambda$ we get
$$
\lim_{t_j\to\infty}t_j\lambda(g)=\lim_{t_j\to \infty}t_j\cF(g,
u)=-\frac{n}{2}.
$$
Using $\lambda\le\int R\,dv/V$, (\Vmono), and the fact that $\lambda\le0$,
we obtain Perelman's monotonicity \cite{P1, \S2} 
$$
\frac{d}{dt}\bar\lambda(g(t))\ge0.
$$
Since $t^{-1}V(t)^{2/n}\to\tilde V_{\infty}^{2/n}$
we obtain
$$
\lim_{t\to \infty} \bar\lambda(g)=-\frac{n}{2}\tilde V_\infty^{n/2}.
$$ 
Our results are summarized in the following Proposition, which by
inspecting the above argument holds even if $\tilde V_\infty=0$.

\proclaim{Proposition \WPboundedprop}
Let $M$ be compact, $g(t)$ a Ricci flow for all $t\ge0$, and $u(t)$
the solution of the conjugate heat equation constructed above.
Then 
$$
\lim_{t\to\infty}\ctW(g(t),u(t),t)=-\log\tilde V_\infty+
\frac{n}{2}(1+\log4\pi),
$$
and
$$
\lim_{t\to\infty}\bar\lambda(g(t))=-\frac{n}{2}\tilde
V^{2/n}_\infty.
$$
\endproclaim

\demo{Remark} 
The upper bound on $\ctW$ can be proven without using
$\cN_+$. If $\tilde V_\infty>0$ then from the monotonicity of
$t^{-n/2}V(t)$ we get the estimate
$$ 
0\le\int_0^\infty t^{-n/2}\int_M\left(R+\frac{n}{2t}\right)\,dv\,
dt<\infty.
$$
So there exist $t^i\to\infty$ so that
$$
(t^i)^{-n/2}\int t^iR\,dv\to-\frac{n}{2}\tilde V_\infty.
$$
Employ this particular sequence $t^i$ in the construction
of $u$ (it should not affect the value of $u$, but we have
not proven this) yields for each fixed $t$
and large $t^i$,
$$
\align 
\ctW(g(t),u^i(t),t)&\le\ctW(g(t^i),u^i(t^i),t^i)\\
&=\frac{1}{V(t^i)}\int t^iR\,dv-\log V(t^i)+\frac{n}{2}\log(4\pi t^i)+n\\
&\to-\log\tilde V_\infty+\frac{n}{2}(1+\log(4\pi)),\qquad t^i\to\infty.\\
\endalign
$$
The result follows by passing $u^i\to u$.
\qeddemo

As a corollary to Proposition \WPboundedprop, 
we obtain the following estimate. 

\proclaim{Proposition \WPestimateprop} Let $M$ be compact and let
$g(t)$, $u(t)$ solve (\RF), (\Ueq) for all $t>0$. Assume $\tilde
V_\infty>0$. Then the rescaled solution $\tilde g(\tilde t):=g(t)/t$,
$\tilde u(\tilde t):=t^{n/2}u$, $\tilde f_+(\tilde t):=f_+(t)$, $\tilde
t=\log t$ satisfies
$$
\int_0^\infty\int_M \tilde u \left|\tilde\Ric+\tilde \covar^2\tilde f_+
+\frac{\tilde g}{2}\right|^2d\tilde v\,d\tilde t<\infty. \tag
\WPestimate
$$
\endproclaim

This estimate may provide enough control to eventually {\it deduce}
convergence of a blowdown sequence in a suitable weak sense. 
We formulate the following natural conjecture as sharply as possible.

\bigskip \proclaim{Conjecure \Hauptvermutung\ (Long-term behavior)}
Let $g(t)$ be a Ricci flow that exists for all $t>0$ on a compact
manifold $M$. Then there is an open set $U\subset M$, 
times $t_j\to\infty$, and a metric $\tilde g_\infty$
such that
$$
g(t_j)/t_j\to\tilde g_\infty\qquad\text{in }U,
$$ 
where $\tilde g_\infty$ is smooth and negative Einstein, the metric
completion of $\tilde g_\infty$ has a singular set of dimension at
most $n-4$, and
$$
\Vol(\tilde g_\infty)=\tilde V_\infty.
$$
\endproclaim

The metric completion of $\tilde g_\infty$ need not be compact or
connected, but has finite volume, so it should be negative Einstein
and not merely an expander. The set $U$ need not be dense in $M$. If
$\tilde V_\infty=0$, then take $U$ to be empty and the statement holds
vacuously; more generally, the assertion is that the scaled volume of the
collapsing parts tends to zero.  The rescaling is designed not to see
these parts: they recede to infinity or hide beyond the singular set.

The conjecture is inspired by the convergence picture in the $3$-dimensional
case \cite{H3, P2}; a good exercise would be to employ the $\ctW$
monotonicity as an alternative way of recognizing the hyperbolic
pieces in that case.  The codimension of the singular set is motivated
by the fact that the right-hand side $\int u|Rm+\cdots|^2$ is unitless
for $n=4$, and is supported by familiar convergence results for Einstein
$4$-manifolds (\cite{A, CC, CT} and many others), as well as analogous
results for mean curvature flow in the critical case \cite{I2} and
other flows.

\bigskip 
\noindent
{\it Formulation of $\mu_+$ and $\nu_+$ functionals.} Define,
analogously to \cite{P1, \S3},
$$
\mu_+(g,\sigma):=\inf_u {\ctW}(g,u,\sigma),\qquad 
\nu_+(g):=\sup_{\sigma>0}\mu_+(g,\sigma),
$$ 
where $u$ varies over functions satisfying $\int_M u\, dv=1$.  We
investigate the existence, smoothness, and monotonicity of these
quantities.  
%As an application, we give another proof of Perelman's
%theorem on expanding breathers.

\proclaim{Theorem \Munuprop} Let $M$ be a closed manifold.

(a) The $\inf$ in the definition
of $\mu_+$ is attained by a unique $u$. 
$\mu_+(g(t),t-T)$ is monotone nondecreasing
under the Ricci flow, and is constant only on an expander with
starting time $t=T$.

(b) If $\lambda(g)<0$, then the $\sup$ in the
definition of $\nu_+$ is attained by a unique $\sigma$. 
$\nu_+(g(t))$ is monotone
nondecreasing  under the Ricci flow, and is constant only on an
expander. 

%(c) (Perelman \cite{P1, \S2}) An expanding breather on $M$ must be
%an expanding soliton (and therefore negative Einstein).
\endproclaim

\demo{Proof}
(a) Write
$$ 
\ctW=\int\left[ \sigma(4|\covar w|^2+Rw^2)+w^2\log w^2\right]\,dv+C,
$$ 
where $w^2=u$. This functional is lower semicontinous and
coercive on $W^{1,2}$, so it possesses a minimizer $w$ subject to
$\int w^2=1$.  Without loss of generality, $w$ is nonnegative. By the
strong maximum principle (taking into account the sign of the
nonlinearity), $w>0$ and therefore $w$ is smooth. The corresponding
$u=w^2$ minimizes the strictly convex functional
$$ 
\ctW=\int\left[ \sigma(4|\covar u^{1/2}|^2+Ru)+u\log u\right]\,dv+C
$$ 
on the cone of smooth functions with $u>0$, $\int u=1$, and
therefore $u$ is unique. It is straightforward to verify that $u$
varies smoothly when $(g,\sigma)$ varies smoothly.  Write
$u_{g,\sigma}$ for this minimizer.

Now let $g(t)$ be a Ricci flow, fix $T\in\R$, and
let $\tilde u(t):=u_{g(t),t-T}$ be the minimizer for each $t$.
Fix $t_0>T$ and let $u(t)$ $t\le t_0$ be the solution of the conjugate
heat equation (\Ueq) with final value $u(t_0)$. 
Since the first variation of $\ctW$ at $t=t_0$ with respect to $u$ 
vanishes, we can compute (smoothly) at $t=t_0$
$$
\frac{d}{dt}\mu_+(g(t),t-T)=\frac{d}{dt}\ctW(g(t),\tilde u(t),t-T)
=\frac{d}{dt}\ctW(g(t),u(t),t-T),
$$
so
$$
\frac{d}{dt}\mu_+(g(t),t-T)
=\int u\left|Rc+\covar^2f_++\frac{g}{2(T-t)}\right|^2\, dv,
$$
where $u$ realizes the minimum at $t=t_0$.
Since $t_0$ was arbitrary, this formula actually holds 
for each $t$, where $u$ realizes the minimum at time $t$. 
If $\mu_+(g(t),t-T)$ is constant on any interval,
the right hand vanishes, so $g(t)$ is an expander with starting time $T$.

(b) By the above, $\mu_+(g,\sigma)$ is continuous in $\sigma$,
$\sigma>0$. We show that it goes to $-\infty$ at the endpoints. Recall
\cite{P1, \S1}
$$
\cF(g,u):=\int \frac{|\covar u|^2}{u}+Ru\,dv,\qquad \lambda(g):=\inf_u\cF(g,w),
$$
where $u$ ranges over functions with $\int u=1$. Fix $u$ and estimate
$$
\mu_+(g,\sigma)\le \ctW(g,u,\sigma)
=\sigma\cF(g,u)+\int u\log u\, dv+\frac{n}2\log(4\pi\sigma)+n.
$$
Thus $\mu_+(g,\sigma)\to-\infty$ as $\sigma\to0$. 
If $\lambda(g)<0$, we select $u$ so that $\cF(g,u)=\lambda(u)$ and
find that $\mu_+(g,\sigma)\to-\infty$ as
$\sigma\to\infty$. In addition it is strictly concave in $\sigma$.
Therefore the supremum in $\nu_+$ is realized by some unique
$\sigma=\sigma_g$. Evidently $\sigma_g$ varies smoothly with $g$. 

Now let $g(t)$ be a Ricci flow and fix $t_0>0$. Since the first
variation of $\mu_+(g(t_0),\sigma)$ is zero at $\sigma=\sigma_{g(t_0)}$,
we may compute (smoothly) at $t=t_0$
$$
\frac{d}{dt}\nu_+(g(t))=\frac{d}{dt}\mu_+(g(t),\sigma_{g(t)})
=\frac{d}{dt}\mu_+(g(t),\sigma+t-t_0).
$$
Thus
$$
\frac{d}{dt}\nu_+(g(t))
=\int u\left|Rc+\covar^2f_++\frac{g}{2\sigma}\right|^2\, dv,
$$
where $(u,\sigma)$ realizes the minimax at $t=t_0$.
Since $t_0$ was arbitrary, this formula actually holds for each $t$,
where $(u,\sigma)$ realizes the minimax at time $t$.
If $\nu_+(g(t))$ is ever constant,
the right hand vanishes, so $g(t)$ is an expander.
\qeddemo

%(c) We prove (c) independent of the above smooth dependence
%on $(g,\sigma)$. Recall that $g(t)$ is called an 
%{\it expanding breather} if there are $t_1<
%t_2$ and a constant $\a>1$ such that $g(t_2)$ 
%differs from $\a g(t_1)$ only by a diffeomorphism. The birth
%time for the breather is 
%$$
%t_0=\frac{\a t_1-t_2}{\a -1}<t_1,
%$$
%which is obtained by solving $(t_2-T)/(t_1-T)=\a$.
%Let $u(x,t_2)$ be the minimizer realizing $\mu_+(g,t_2-T)$,
%and let $u(x,t)$, $t\le t_2$, solve the
%%conjugate heat equation with ending data
%$u(x,t_2)$. Using the monotonicity of $\ctW$, the scaling
%$\mu_+(kg,k\sigma)=\mu_+(g,\sigma)$, and the fact that $g(t)$ is a breather,
%we obtain
%$$
%\split
%\mu_+(g(t_1),t_1-T)&\le {\ctW}(g(t_1),u(t_1),t_1-T)\\
%&\le {\ctW}(g(t_2,u(t_2),t_2-T)\\
%&=\mu_+(g(t_2), t_2-T)\\
%&=\mu_+(\alpha g(t_1), \alpha(t_1-T))\\
%&=\mu_+(g(t_1),t_1-T).\\
%\endsplit
%$$
%Thus the right-hand side of the monotonicty formula vanishes,
%so $g(t)$ is an expander.
%\qeddemo

\bigskip
\subheading{\S2 Forward reduced distance and reduced volume}

\bigskip

In this section we derive the dual version of Perelman's
monotonicity on the reduced volume in \cite{P1, \S7}.

Let $g(t)$ solve the Ricci flow on $M\times [0,
T]$. Fix $x_0$ and let $\gamma$ be a path $(x(\eta),\eta)$
joining $(x_0, 0)$ to $(y, t)$. 
Analogous to \cite{P1}, we define
$$ 
\ctL(\gamma)=\int_0^t\sqrt{\eta}
\left({\Rscalar}+|\gamma'(\eta)|^2\right)\, d\eta.\tag 2.1
$$
Let $X\equiv\gamma'(t)$ 
and let $Y$ be a variational vector field along $\gamma$. 
One calculates the first variation of $\ctL$ to be
$$
\delta {\ctL} =\left.2\sqrt{t}\langle X, Y\rangle\right|^{t}_{0}
+\int_0^t\sqrt{\eta}\langle Y,
\covar  R-2\covar_X X +4\Ric(X, \cdot)-\frac{1}{\eta}X\rangle\, d\eta.\tag
2.2
$$
From
this one can write the ${\ctL}$-geodesic equation:
$$
\covar_X X -\frac{1}{2}R +\frac{1}{2t}X -2\Ric(X, \cdot)=0. \tag
2.3
$$
This differs from \cite{P1} only by the sign in front of $\Ric(X,
\cdot)$. Let $L_+(y, t)$ denote the length of a shortest 
$\ctL$-geodesic joining $(x_0, 0)$ to $(y, t)$. 
From (2.2) we know that
$$
\covar  L_+ =2\sqrt{t}X(t). \tag 2.4
$$
Following closely the computation of  \cite{P1, \S7} we obtain
$$
\gather
|\covar L_+|^2=-4t{\Rscalar} +4t({\Rscalar}+|X|^2) \tag 2.5\\
\frac{\pa L_+}{\pa t}=2\sqrt{t}{\Rscalar} -\sqrt{t}({\Rscalar}+|X|^2). \tag 2.6
\endgather
$$
A calculation using the geodesic equation shows that
$$
\frac{d}{d t} \left({\Rscalar}+|X|^2\right)=H(X)-\frac{1}{t}
\left({\Rscalar}+|X|^2\right),
$$
where $H(X):=\pa\Rscalar/\pa t+2\langle\covar {\Rscalar},
X\rangle+2\Ric(X,X)+\Rscalar/t$  is the exactly twice
the traced Li-Yau-Hamilton differential Harnack expression in
\cite{H1}. (Notice that the $H$ in \cite{P1, \S7} 
also equals the LYH expression, but evaluated at a {\it negative}
time $t=-\tau$.) 
This gives
$$
t^{3/2}\left({\Rscalar}+|X|^2\right)=K +\frac{1}{2}L_+ \tag 2.7
$$
where
$$
K:=\int_0^t \eta^{3/2}H(X)\, d\eta.
$$
We then have
$$
\gather
|\covar L_+|^2=-4t{\Rscalar}+\frac{2}{\sqrt{t}}L_+ +\frac{4}{\sqrt{t}}K,
\tag 2.8 \\
\frac{\pa L_+}{\pa t}=2\sqrt{t}{\Rscalar}-\frac{1}{t}K -\frac{1}{2t} L_+. \tag 2.9
\endgather
$$ 
These differ from \cite{P1} only by the sign in front of $K$ (but
recall that the interpretation of $H$ also differs).

We can similarly do the second variation computation. For the
a variation vector $Y$ satisfying $Y(0)=0$, we get
$$
\split
\delta^2_Y\ctL&=2\sqrt{t}\langle\covar_Y Y, X\rangle+\int_0^t
\sqrt{\eta}\left(Y\cdot Y\cdot {\Rscalar}
-2\langle\covar_Y Y,\covar_XX\rangle\right.\\
&\quad\left.+4Y(\Ric(Y, X))-2X(\Ric(Y, Y))-2R(X, Y, X, Y)\right.\\
&\quad\left.+2|\covar_XY|^2
-\frac{1}{\eta}\langle\covar_Y Y, X\rangle\right)d\eta.
\endsplit
\tag 2.10
$$
Using the ${\ctL}$-geodesic equation we get
$$
\split \delta^2_Y\ctL&=2\sqrt{t}\langle\covar_Y Y, X\rangle+\int_0^t
\sqrt{\eta}\left(\covar^2 {\Rscalar}(Y, Y)-2R(X, Y, X, Y)\right.\\
&\quad \left.+2|\covar_X Y|^2 4\covar_Y \Ric( Y, X)
-2\covar_X \Ric(Y, Y)\right)\, d\eta.
\endsplit
\tag 2.11
$$
We choose $Y$ such that
$$
\covar_X Y =\Ric(Y, \cdot)+\frac{1}{2t}Y, \tag 2.12
$$
which in particular implies that  $|Y(\eta)|^2=\eta/t$.
Then we obtain 
$$
\split \covar^2L_+(Y, Y)&\le \frac{1}{\sqrt{t}}+\int_0^t
\sqrt{\eta}\left(\covar^2 {\Rscalar} (Y, Y)-2R(X, Y, X,
Y)+2|\Ric(\cdot,Y)|^2+\frac{2}{\eta}\Ric(Y, Y)\right.\\
&\quad \left. (4\covar_Y \Ric(Y, X)-4\covar_X \Ric(Y, Y))+2\covar_X
\Ric(Y, Y)\right)\, d\eta. \endsplit \tag 2.13
$$
Now we use the fact that
$$
\frac{d}{d t}\left(\Ric(Y, Y)\right)=\frac{\pa \Ric}{\pa t}(Y, Y)
+\covar_X \Ric(Y,
Y)+2|\Ric(Y, Y)|^2+\frac{1}{t}\Ric(Y, Y)
$$
which implies 
$$
\split 2\int_0^t \sqrt{\eta}&\covar_X \Ric(Y, Y)\, d\eta=2\sqrt{t}\Ric(Y, Y)\\
&-\int_0^t \sqrt{\eta}\left(3\frac{\Ric(Y,
Y)}{\eta}-2\frac{\pa \Ric}{\pa t}(Y, Y)-4|\Ric(Y, \cdot)|^2\right)\, d \eta.
\endsplit\tag 2.14
$$
Plugging (2.14) into (2.13) we have that
$$
\covar^2L_+(Y, Y)\le \frac{|Y|^2}{\sqrt{t}}+2\sqrt{t}Rc(Y,Y)-\int_0^t
\sqrt{\eta} H(X, Y)\, d\eta. \tag 2.15
$$
Here
$$
\split
H(X, Y):=&-\covar^2 {\Rscalar}(Y, Y)+2R(X, Y, X,Y)+2|\Ric(\cdot, Y)|^2\\
&\qquad+\frac{1}{t}\Ric(Y, Y)+2\frac{\pa \Ric}{\pa t}(Y,Y)\\
&\qquad-4\covar_Y \Ric(Y, X)+4\covar_X \Ric(Y, Y).
\endsplit
\tag 2.16
$$
This is exactly twice Hamilton's matrix LYH expression. We get the following.

\proclaim{Corollary \PosCurvOp} If $(M, g(t))$ has nonnegative curvature
operator, then
$$
\covar^2L_+(Y, Y)\le \frac{|Y|^2}{\sqrt{t}}+2\sqrt{t}Rc(Y,Y)
$$
\endproclaim

%note: we don't know any applications of this

\noindent
In the general curvature case, tracing (2.15) in $Y$ yields
$$
\D L_+ \le \frac{n}{\sqrt{t}}+2\sqrt{t}{\Rscalar}-\frac{1}{t}K. \tag
2.17
$$
Now defining
$$
\ell_+(y, t):=\frac{1}{2\sqrt{t}}L_+(y, t),\qquad \bar L_+(y,
t):=4t\ell_+(y, t),
$$
we get
$$
\gather
|\covar \ell_+|^2=-{\Rscalar}+\frac{\ell_+}{t}+\frac{K}{t^{3/2}}, \tag
2.8$^\prime$\\
\frac{\pa \ell_+}{\pa t}
={\Rscalar} -\frac{K}{2t^{3/2}}-\frac{\ell_+}{t}, 
\tag 2.9$^\prime$\\
\D \ell_+\le {\Rscalar}+\frac{n}{2t} -\frac{K}{2t^{3/2}}. \tag
2.17$^\prime$
\endgather
$$
These then imply
$$
\gather
\frac{\pa \ell_+}{\pa t} +\D \ell_+ 
+|\covar \ell_+|^2-{\Rscalar}-\frac{n}{2t}\le 0, \tag 2.18\\
\heat  \left(\bar L_++2nt\right)\ge 0, \tag 2.19
\endgather
$$
and negativity of a quantity resembling $-v_+$:
$$
t(2\D \ell_++|\covar \ell_+|^2-{\Rscalar})-\ell_+-n\le 0. \tag 2.20
$$ 
Note that $\ell_+$ resembles
$-f_+$, in contrast to the shrinker case, where $\ell$ is like $f$
\cite{P1, \S9, CHI}. We have the following theorem from (2.18), (2.19).

\proclaim{Theorem \ExpanderEll} Let $g(t)$ solve the Ricci
flow on $M\times[0, T]$. Then
$$
\hat u(x,t):=\frac{e^{\ell_+(x,t)}}{(4\pi t)^{n/2}}
$$
is a super-solution to the conjugate heat equation, that is,
$
\p\hat u/\p t+\D\hat u-{\Rscalar}\hat u\le 0$.
In particular, when $M$ is compact, the forward reduced volume
$$
\theta_+(t)=\theta^{(x_0,0)}_+(t):= \int_M \hat u\, dv \tag 2.21
$$
is monotone non-increasing along the flow.
Furthermore, $\bar L_++2nt$ is a super-solution to the standard heat equation. 
\endproclaim

Remarkably, the entropies $\cW$ and $\ctW$ both increase in forward
time, whereas the reduced volumes increase as $t$
approaches the reference time $t=0$. 

The forward reduced volume is a constant precisely when
$g(t)$ is a compact expanding soliton and $(x_0,0)$ is the ``vertex'' of the
soliton. In this case we get
$$
\log\theta_+(t)=-\nu_+(g),
$$
a constant independent of time. Details will be given in future work;
see \cite{CHI} for the shrinker case.
The following proposition says that $\theta_+$ is bounded away from zero
when $\tilde V_\infty>0$, in the same way as the other monotone quantities.

\proclaim{Proposition \ThetaLower} The forward reduced volume is
bounded below by
$$
\log\theta(t)\ge \log \tilde V(t)-\frac{n}{2}(1+\log 4\pi).
$$
\endproclaim

\demo{Proof}
Since $R\ge -n/2t$, we get for any curve $\gamma$ from $(x_0,0)$ to $(p,t)$,
$$
\ctL(\gamma)\ge\int_0^t\sqrt{\eta}\left(-\frac{n}{2\eta}\right)\,d\eta
=-n\sqrt{t},
$$
so $\ell_+(p,t)\ge -n/2$ and 
$$
\theta_+(t)\ge \frac{V(t)e^{-n/2}}{(4\pi t)^{n/2}}
=\frac{\tilde V(t)}{(4\pi e)^{n/2}}. 
\qed
$$
\enddemo

This also shows that $\theta_+(t)\to\infty$ as $t\to0$ unless
$V(t)\sim t^{n/2}$ for small $t$ (maximal volume decay). 
In particular, when the initial metric is
smooth except for a conic singularity $x_0$, as may arise
by evolution from a smooth metric defined for negative times \cite{FIK},
then $\theta_+(t)$ will not be
useful in proving that a blowup about the vertex $(x_0,0)$
for small $t>0$ yields an expander.

\end